\documentclass[11pt, final]{amsart}
\usepackage[cp850]{inputenc}
\usepackage{amssymb}
\usepackage[mathscr]{eucal}

\theoremstyle{plain}
\newtheorem{theorem}{Theorem}
\newtheorem{lemma}{Lemma}

\newtheorem{example}{Example}

\theoremstyle{remark}
\newtheorem{remark}{Remark}

\def\cA{\mathscr{A}}
\def\cS{\mathscr{S}}
\def\Ref{\operatorname{ref}}
\def\Iso{\operatorname{Iso}}
\def\Aut{\operatorname{Aut}}

\title{Reflexivity of the isometry group of some classical spaces}
\author{F. Cabello S\'{a}nchez}
\address{Departamento de Matem\'{a}ticas, Universidad de Extremadura,
Avenida de Elvas, 06071-Badajoz, Spain}
\email{fcabello@unex.es}
\thanks{The first author was supported in part by DGICYT project PB97-0377
and HI project 1997-0016}
\author{L. Moln\'{a}r}
\address{Institute of Mathematics and Informatics,
Lajos Kossuth University, H-4010 Debrecen, P.O. Box 12, Hungary}
\email{molnarl@math.klte.hu}
\thanks{The second author was supported by the Hungarian National
Foundation
for Scientific Research (OTKA) and by the Ministry of Education (FKFP)}
\thanks{1991 Mathematics Subject Classification 47B49, 46B04, 47B10.}
\begin{document}
\begin{abstract}
We investigate the reflexivity of the isometry group and the automorphism
group of some important metric linear spaces and
algebras. The paper consists of the following sections:
1. Preliminaries. 2.
Sequence spaces. 3. Spaces of measurable functions. 4. Hardy spaces.
5. Banach algebras of holomorphic functions. 6. Fr\'echet\break algebras of
holomorphic functions. 7. Spaces of continuous functions. \end{abstract}
\maketitle
\bigskip\bigskip

\begin{center}{\scshape Introduction}
\end{center}
This paper is concerned with the reflexivity of the isometry group and the
automorphism group of certain particular but important topological vector
spaces and algebras. Although we deal mainly with Banach spaces, we are
also
interested in other (not necessarily locally convex) metric linear spaces
and some Fr\'{e}chet algebras.

Reflexivity problems for subalgebras of the algebra of all bounded
linear operators acting on a Hilbert space represent one of the most
active research areas in operator theory. The study of similar questions
concerning sets of linear transformations on Banach algebras rather than on
Hilbert spaces was initiated by Kadison \cite{Kad} and Larson \cite{Lar}.
In \cite{Kad}, motivated by the study of Hochschild cohomology of operator
algebras, the reflexivity of the Lie algebra of all derivations on a von
Neumann algebra was treated. In \cite[Some concluding remarks (5), p.
298]{Lar}, Larson raised the question of the reflexivity of the
automorphism
group of Banach algebras. This problem
was investigated for several algebras in \cite{bm, Mol1, Mol2, MGy,
mz2}. The present article is a continuation of that work.

We describe the results of the paper as follows.
The first section is preliminary.
In Section \ref{ssequence} we investigate symmetric spaces. We prove that,
with the sole exception of $l_2$, every $F$-space with a symmetric
basis has algebraically reflexive isometry group. Curiously enough, the
isometry group of any non-separable ``symmetric" space fails to be
algebraically reflexive. Section \ref{slebesgue} concentrates on Lebesgue
spaces. We show the extreme nonreflexivity of the isometry group of the
spaces $L_p(\mu)$ $(0<p<\infty)$
for homogeneous measures $\mu$, thus obtaining that the only infinite
dimensional
Lebesgue spaces whose isometry groups are reflexive are the sequence spaces
$l_p$ for $p\neq2$. In contrast, the isometry groups of the Hardy spaces
$H^p$ $(0<p<\infty)$ are topologically reflexive for all $p\neq2$. This
will
be proved in Section \ref{shardy}.

Sections \ref{sbounded} and \ref{sunbounded} deal with algebras of
holomorphic functions. We prove that the isometry group and the
automorphism
group of the disc algebra are topologically reflexive. The same is true for
$H^\infty(\Omega)$, $\Omega$ being any simply connected domain in the
plane.
Furthermore, we consider the "unbounded" case: it is shown that the
automorphism
group of the Fr\'{e}chet algebra $H(\Omega)$ is topologically reflexive if
and only if $\Omega\neq\mathbb C$.

Finally, we study algebras of continuous functions.
We solve some problems posed in \cite{mz} by presenting Banach
spaces whose isometry
groups are either trivial or very large. We give an example of a compact
Hausdorff
space $K$ such that the isometry group and the automorphism group
of the Banach algebra $C(K)$ both fail to be reflexive, thus
verifying a conjecture formulated in \cite{mz}, where it was proved that
the
isometry group and the automorphism group of $C(K)$ are algebraically
reflexive in case $K$ is first countable. Furthermore, we
exhibit a Banach space $X$ with the property that its isometry group is
topologically
reflexive but the isometry group of its dual space $X^*$ is even
algebraically nonreflexive. Next, we present another Banach space $Y$ whose
isometry group
is not algebraically reflexive and yet the isometry group of $Y^*$
is algebraically reflexive.

\section{Preliminaries}

Let $X$ be a topological vector space and let $B(X)$ be the algebra of all
continuous linear operators on $X$. Given any subset $\cS\subset B(X)$,
define
$$
\Ref_{al} \cS=\{T\in B(X): Tx\in \cS x\text{ for all }x\in X\}$$ $$
\Ref_{to} \cS=\{T\in B(X): Tx\in{\overline {\cS x}}\text{ for all }x\in
X\},
$$
where $\cS x= \{Lx:L\in \cS\}$ and the bar stands for the closure in $X$.
The
set $\cS$ is said to be algebraically reflexive if $\Ref_{al}\cS= \cS$ and,
similarly, $\cS$ is called topologically reflexive if $\Ref_{to} \cS= \cS$.
Thus, reflexive sets of operators are, in some sense, completely determined
by their local actions on the underlying space. Sometimes the operators in
$\Ref_{al} \cS$ are said to belong locally to $\cS$.\medskip

Since there is no clear intrinsic reason to restrict our attention to the
locally convex setting when dealing with local surjective isometries, we
consider $F$-spaces, not only Banach spaces.

Recall from
\cite[p. 2]{kpr} that a $\Delta$-norm on a real or complex vector space $X$
is a non-negative real-valued function on $X$ satisfying \begin{itemize}
\item[(1)]\quad $\|x\|>0$ for all $0\neq x\in X$.
\item[(2)]\quad $\|\alpha x\|\leq \|x\|$ for all $x\in X$ and all
$\alpha\in\mathbb K$ with $|\alpha|\leq1$.
\item[(3)]\quad $\lim_{\alpha\to
0}\|\alpha x\|=0$ for all $x\in X$.
\item[(4)]\quad $\|x+y\|\leq
K(\|x\|+\|y\|)$ for some constant $K$ independent of $x,y\in X$.
\end{itemize}A $\Delta$-norm on $X$ induces a metrizable linear topology
for which the sets $U_n=\{x\in X:\|x\| <1/n\}$ form a neighbourhood base
at the origin and, conversely, every linear metrizable topology comes
from a $\Delta$-norm. An $F$-norm is a $\Delta$-norm satisfying
\begin{itemize}
\item[(5)]\quad$\|x+y\|\leq\|x\|+\|y\|$ for all $x,y\in X$. \end{itemize}
Any $F$-norm induces a translation-invariant metric on $X$ in the obvious
way
and every invariant metric compatible with the linear structure is induced
by some $F$-norm. An $F$-space is a complete $F$-normed space.

Finally, a quasi-norm is a $\Delta$-norm which is homogeneous in the sense
that
\begin{itemize}
\item[(6)]\quad$\|\alpha x\|=|\alpha|\|x\|$ for all $x\in X, \alpha\in
\mathbb K$.
\end{itemize}
Observe that (6) implies both (2) and (3) so that quasi-norms can be
defined
by (1), (6) and (4). A quasi-normed space is a vector space together with
some specified quasi-norm. Such a space is locally bounded, that is, it has
a bounded neighbourhood of zero. Conversely, every locally bounded
topology is induced by a quasi-norm. A quasi-Banach space is a complete
quasi-normed space.

We denote by $\Iso(X)$ the group of all surjective (linear) isometries of
the $\Delta$-normed space $X$. Also, when $\cA$ is a topological algebra,
$\Aut(\cA)$ denotes the group of all continuous automorphisms of $\cA$. In
accordance with what is written above, we call the elements of
$\Ref_{al}(\Iso(X))$ and $\Ref_{al}(\Aut(\cA))$ local surjective isometries
and local automorphisms, respectively.

One little problem with $\Delta$-norms is that a $\Delta$-norm need not be
continuous with respect to the topology it induces. (The continuity of
a norm is a consequence of the triangle inequality.) This has some
unpleasant consequences. For instance, the operators in
$\Ref_{to}(\Iso(X))$
need not be into isometries. An interesting class of continuous quasi-norms
is that of the so-called $p$-norms $(0<p\leq1)$. These are quasi-norms
satisfying
$$
\|x+y\|^p\leq\|x\|^p+\|y\|^p,
$$
from which continuity immediately follows. Clearly,
if $\|\cdot\|$ is a $p$-norm, then $\|\cdot\|^p$ is an $F$-norm with the
same isometries.

\section{Sequence spaces}\label{ssequence}
In this section we study sequence spaces. A basis of an $F$-space $X$ is a
sequence $(e_n)$ so that every $x\in X$ has a unique expansion $
x=\sum_{n=1}^\infty x_ne_n.
$
A basis $(e_n)$ is said to be symmetric if $$
\left\|\sum_{n=1}^\infty
x_ne_n\right\|=\left\|\sum_{n=1}^\infty\varepsilon_nx_ne_{\pi(n)}\right\|
$$
holds for every choice of scalars $\varepsilon_n$ of modulus 1 and every
permutation $\pi$
of the positive integers. Perhaps the most interesting class of spaces with
symmetric basis is that of the Orlicz sequence spaces
(see \cite{kpr, lt, ro} for definitions) which contains the Banach spaces
$l_p$ for $1\leq p<\infty$, the quasi-Banach spaces $l_p$ for $0<p<1$ and
other locally bounded and even non-locally bounded spaces. Another
important
class related to $l_p$ spaces is that of the Lorentz sequence spaces
\cite{lt}.

\begin{theorem}\label{symbasis} Let X be an F-space not isomorphic to
$l_2$.
Suppose that X has a symmetric basis. Then the isometry group of X is
algebraically reflexive. \end{theorem}

\begin{proof} By a result of Rolewicz \cite[Theorems 9.8.3 and 9.8.5]{ro}
(see also \cite{gl}), every surjective isometry of $X$ is of the form $$
\sum x_ne_n\longmapsto\sum \varepsilon_nx_ne_{\pi(n)}, $$
where $|\varepsilon_n|=1$ and $\pi$ is a permutation of $\mathbb N$.
It follows that if $T$ is a local surjective isometry
of $X$, then there is an injective mapping $\varphi$
on $\mathbb N$ for which $$
Te_n=\sigma_ne_{\varphi(n)},
$$
where $|\sigma_n|=1$. Hence $T$ is given by
$$
T\left(\sum x_ne_n\right)=\sum\sigma_nx_ne_{\varphi(n)}. $$
The theorem will be proved if we show that $\varphi$ is surjective. It is
easily seen that $X$ contains an $x=\sum x_ne_n$ with $x_n>0$, $x_n\neq
x_m$
$(n\neq m)$. Taking
$L\in\Iso(X)$ so that $Tx=Lx$, one obtains that $$
\sum\sigma_nx_ne_{\varphi(n)}=\sum\varepsilon_nx_ne_{\pi(n)} $$
for some permutation $\pi$ of $\mathbb N$. Hence $\varphi=\pi$, $\varphi$
is
surjective and the proof is complete. \end{proof}

\begin{remark}
Observe that, in general, we do not have topological reflexivity in
Theorem~\ref{symbasis}. In fact, let $X=c_0$ and consider the unilateral
shift $S$ on $X$. It is easy to see that for every $x\in X$ and $\epsilon
>0$ there exists a surjective isometry $L$ such that $\|Sx-Lx\| <\epsilon$.
Therefore, $S\in \Ref_{to}(\Iso(X))$ but $S$ is not surjective.
\end{remark}

The following example shows that separability is essential in
Theorem \ref{symbasis}.

\begin{example}\label{uncountable} Let X be a $\Delta$-normed space of
functions on an index set $\Gamma$. Suppose that \begin{itemize}
\item[(a)] for every $f\in X$ the
set $\{\gamma\in\Gamma: f(\gamma)\neq0\}$ is at most countable,
\item[(b)] for every bijection $\varphi$ of $\Gamma$, the map $f \mapsto
f\circ \varphi$ is a surjective isometry of X.
\end{itemize}
If $\Gamma$ is
uncountable, then the isometry group of X is algebraically nonreflexive.
\end{example}

\begin{proof} We closely follow \cite[last Remark]{mz2}. Let $\Lambda$
be
a proper subset of $\Gamma$ with a bijection $\varphi:\Lambda\to\Gamma$.
Define $T:X\to X$ by
$$
(T f)(\gamma)=\begin{cases}f(\varphi(\gamma))\quad&\text{if $\gamma$
belongs to $\Lambda$,}\\
0&\text{elsewhere.}
\end{cases}
$$
By (b), $T$ is obviously non-surjective. To see that
$T\in\Ref_{al}(\Iso(X))$,
fix $f\in X$. Since $f$ has at most countable support one can find a
bijection
$\phi:\Gamma\setminus\varphi^{-1}(\text{supp}(f))\to\Gamma\setminus
\text{supp}(f)$. Define a bijection on $\Gamma$ by $$
\varphi_f(\gamma)=\begin{cases}\varphi(\gamma)\quad&\text{if }
\gamma\in\varphi^{-1}(\text{supp}(f)) \\ \phi(\gamma)&\text{elsewhere.}
\end{cases}
$$
Let $T_f$ be given on $X$ by $T_f(g)=g\circ\varphi_f$. Then
$T_f\in\Iso(X)$ and $T_f(f)=T(f)$. This completes the proof.
\end{proof}

\begin{remark} It is well-known that the isometry group of any infinite
dimensional real or complex Hilbert space fails to be algebraically
reflexive. This is because Hilbert spaces are isotropic: given $x,y\in H$
with $\|x\|=\|y\|$ there is $T\in \Iso(H)$ such that $y=Tx$. Hence
$\Ref_{al}(\Iso(H))$ is as large as possible and
contains all into isometries. Therefore, $\Iso(H)$ cannot be reflexive
unless $H$ is finite dimensional. In fact, the isometry group of any
infinite dimensional complex Hilbert space is algebraically nonreflexive
not
only with respect to the original Hilbert space norm, but also
with respect to the so-called spin norms \cite[Theorem 3.7]{mz}. It is
natural to ask if the same occurs with any equivalent norm (or renorming,
in
short). Although in \cite{mz} an affirmative answer was conjectured, the
following result shows that the answer is strongly negative. \end{remark}

\begin{theorem} Every Banach space admits a renorming whose isometry group
is topologically reflexive.
\end{theorem}

\begin{proof} Let $Y$ be a Banach space. By a result of Jarosz \cite{Ja},
there is a renorming $X$ of $Y$ with trivial isometries, that is, such that
$\Iso(X)=\{\sigma I : \sigma\in\mathbb{K}, |\sigma|=1\}$. Clearly,
$\Iso(X)$
is topologically reflexive.
\end{proof}

\begin{remark}
Suppose that $Y$ is a Hilbert space. Using \cite{St} instead of \cite{Ja}
one obtains that for every $\epsilon>0$ there is a Banach space $X$ with
trivial isometries and $\epsilon$-isometric to $Y$ in the sense that there
is an isomorphism $T:X\to Y$ with
 $\|T\|\cdot\|T^{-1}\|<1+\epsilon$. This means that the new norm can be
chosen to be a very small perturbation of the original Hilbert space norm
of $Y$.
\end{remark}

\section{Spaces of measurable functions}\label{slebesgue} Let
$(\Omega,\Sigma,\mu)$ be a
measure space. For $0<p<\infty$, define $L_p(\mu)$ to be the space of all
real or complex measurable functions on $\Omega$ for which $$
\|f\|_p=\biggl(\int_\Omega|f|^pd\mu\biggr)^{1/p} $$
is finite with
the usual convention about identifying functions equal almost everywhere.
Observe that $\|\cdot\|_p$ is a norm only if $p\geq 1$. For $0<p<1$ it is
only a quasi-norm (in fact, a $p$-norm) and $L_p(\mu)$ is a quasi-Banach
space.

By a famous theorem of Maharam, every Lebesgue space (that is, $L_p$ space)
is isometrically representable
as
$$
L_p(\mu)=l_p(\Gamma)\oplus_p\biggr(\sum_{i\in
I}L_p(\lambda^{{\mathfrak{c}}_i})\biggl)_p $$
where $\Gamma$ and $I$ are
(possibly empty) sets, ${\mathfrak{c}}_i$ are infinite cardinals and
$\lambda$ denotes the Lebesgue measure on $[0,1]$ (observe that, for
instance,
$L_p(\lambda)=L_p(\lambda^{\omega}))$. The subscript $p$ indicates that the
corresponding direct sum is taken in the $l_p$ sense.

Clearly, $\mu$ is $\sigma$-finite if and only if $\Gamma$ and $I$ are
countable sets. Also, $\mu$ is homogeneous in the sense of \cite{gjk} if
and
only if $\Gamma$ is empty and all ${\mathfrak{c}}_i$ coincide (apart from
the trivial case when $I$ is empty and $\Gamma$ is a singleton). Homogenity
means that $\Sigma|_{\Omega'}$ and $\Sigma|_{\Omega''}$ are Boolean
isomorphic whenever $\Omega'$ and $\Omega''$ are subsets of $\Omega$ with
positive finite measure, which implies that, for every $0<p<\infty$, the
spaces $L_p(\Omega',\mu)$ and $L_p(\Omega'',\mu)$ are isometrically
isomorphic. In that case the isometry group of $L_p(\mu)$ has at most two
orbits on the unit sphere. More precisely, given $f,g\in L_p(\mu)$ with
$\|f\|_p=\|g\|_p\neq 0$, there is a surjective isometry of $L_p(\mu)$
mapping $f$ into $g$ if and only if either both $f$ and $g$ are nonzero
almost everywhere or both $f$ and $g$ vanish on sets of positive measure
(this follows from \cite[Lemma 1.4]{gjk}).
This obviously implies that these spaces are almost isotropic: given
$f,g\in
L_p(\mu)$ with $\|f\|_p=\|g\|_p=1$ and $\varepsilon>0$ there is $T\in\Iso(
L_p(\mu))$ fulfilling $\|g-Tf\|_p\leq\varepsilon$. If, in addition to be
homogeneous, $\mu$ is not $\sigma$-finite, then $L_p(\mu)$ is isotropic for
all $0<p<\infty$.

\begin{theorem}\label{lebesgue}
Let X be an infinite dimensional Lebesgue space.
Then $\Iso(X)$ is algebraically reflexive if and only if $X=l_p(\mathbb N)$
with $p\neq2$. \end{theorem}

\begin{proof}
The standard basis of $l_p(\mathbb N)$ is symmetric for every
$0<p<\infty$, so Theorem
\ref{symbasis} implies that $\Iso(l_p(\mathbb N))$ is algebraically
reflexive for every $p\neq2$.

For the converse we need the following lemma whose easy proof is left to
the reader.

\begin{lemma}
Let X be a $\Delta$-normed space whose isometry group is algebraically
reflexive and let Y be a linear subspace of X. Suppose that Y has
a complement Z in X such that $\|y+z\|=\Phi(\|y\|,\|z\|)$ holds for some
function $\Phi:\mathbb R^2\to\mathbb R$ and all $y\in Y, z\in Z$. Then
$\Iso(Y)$ is algebraically reflexive too. \end{lemma}

Now, it clearly suffices to see that $\Iso(X)$ is algebraically
nonreflexive for $X$ equal either $L_p(\lambda^{\mathfrak c})$ or
$l_p(\Gamma)$
with $\Gamma$ uncountable and then apply Lemma 1. But $l_p(\Gamma)$ is just
a particular case of Example \ref{uncountable}. The following result ends
the proof of Theorem~\ref{lebesgue}. \end{proof}

\begin{lemma}\label{homogeneous} For any cardinal ${\mathfrak c}$ the
isometry group of the space $L_p(\lambda^{\mathfrak c})$ is algebraically
nonreflexive.
\end{lemma}

\begin{proof} First observe that two functions $f$ and $g$ in an arbitrary
$L_p$
space with $p\neq2$ have almost disjoint supports (that is, $fg=0$ holds
almost everywhere) if and only if
$$
\|f+g\|_p^p+\|f-g\|_p^p=2(\|f\|_p^p+\|g\|_p^p). $$
It follows that if $T$ is an into isometry between $L_p$ spaces and $f$
vanishes on a set of positive measure, then so does $Tf$. Hence, if
in addition $T$ is
surjective, then $Tf$ is nonzero almost everywhere if and only if $f$ is
nonzero almost everywhere.

Thus, in view of the structure of the orbits of the isometry group of
$L_p(\mu)$ for a homogeneous measure $\mu$, it is clear that an into
isometry of $L_p(\lambda^{\mathfrak c})$ is locally surjective if and only
if it preserves the almost everywhere nonzero functions as well as the
functions with support of positive measure.

Let $\Omega=[0,1]^{\mathfrak c}$. Write
$\Omega=\Omega'\oplus\Omega''$ with $\Omega', \Omega''\in\Sigma$ of
positive measure. Obviously, $L_p(\Omega,\lambda^{\mathfrak
c})=L_p(\Omega',\lambda^{\mathfrak
c})\oplus_pL_p(\Omega'',\lambda^{\mathfrak c})$. By the homogenity of
$\lambda^{\mathfrak c}$, there are surjective isometries
$T':L_p(\Omega,\lambda^{\mathfrak c})\to L_p(\Omega',\lambda^{\mathfrak
c})$
and $T'':L_p(\Omega,\lambda^{\mathfrak c})\to
L_p(\Omega'',\lambda^{\mathfrak c})$.
Define $T:L_p(\lambda^{\mathfrak c})\to L_p(\lambda^{\mathfrak c})$ as $$
T(f)=\frac{T'(f)+T''(f)}{2^{1/p}}.
$$
Clearly, $T$ is a non-surjective into isometry. On the other hand, it is
easily seen that $Tf$ is nonzero almost everywhere if and only if $f$
is. Therefore, $T$ is a local surjective isometry of
$L_p(\Omega,\lambda^{\mathfrak c})$. This proves the lemma. \end{proof}

\begin{remark} Let $\varphi$ be an Orlicz function (see \cite[p. 29]{kpr}).
Then $L_\varphi(\mu)$, the Orlicz function space determined by $\varphi$ on
the finite measure space $(\Omega,\Sigma,\mu)$ is the space of all
measurable functions on $\Omega$ for which $$
\|f\|_\varphi = \int_\Omega\varphi(|f(\omega)|)d\mu $$
is finite. It can be proved that $\|\cdot\|_\varphi$ is a complete
$\Delta$-norm on $L_\varphi(\mu)$. Under some additional hypotheses (e.g.
concavity) $\|\cdot\|_\varphi$ is even an $F$-norm. It would be interesting
to know whether the isometry group of the spaces $L_\varphi(\mu)$ is
reflexive. For instance, we do not known whether the isometry group of the
space of all measurable functions on $[0,1]$ is algebraically reflexive
(that space is determined by the concave Orlicz function
$\varphi(t)=t/(1+t)$). It should be noted that, with the exception of the
Lebesgue spaces $L_p(\mu)$ (which correspond to the Orlicz functions
$t\mapsto t^p$), every surjective isometry on a "reasonable"
$L_\varphi(\mu)$ is induced by a measure-preserving automorphism of
$\Sigma$
(see \cite{lam}).
\end{remark}

\section{Hardy spaces}\label{shardy}
In the sequel, we denote by $U$ the open unit disc in the plane and by
$\mathbb T$ the unit circle. For $0<p<\infty$, the Hardy space $H^p$ is the
space of all holomorphic functions $f:U\to\mathbb C$ for which
$$
\|f\|_p=\sup_{0\leq
r<1}\bigg(\frac{1}{2\pi}\int_0^{2\pi}|f(re^{i\theta})|^pd\theta\bigg)^{1/p}
$$
is finite (the space $H^\infty$ will be treated in the next section). If
$p\geq1$, then $\|\cdot\|_p$ is a complete norm on $H^p$, while for $0<p<1$
it is only a $p$-norm and $H^p$ is a quasi-Banach space.

We refer the reader to \cite{du} for general information about Hardy
spaces. Here we only recall the well-known inequality
$$
|f(re^{i\theta})|\leq2^{1/p}(1-r)^{-p}\|f\|_p\quad\quad(r<1, \:f\in H^p),
$$
(see \cite[p. 36, Lemma]{du} which implies that for every $z\in U$, the
point evaluation $f\mapsto f(z)$ is continuous on $H^p$.

\begin{theorem}For every $p\neq2$ the isometry group of $H^p$ is
topologically reflexive.
\end{theorem}

Before going into the proof, recall from \cite{fo} that every into isometry
of $H^p$ $(p\neq2)$ has the form
\begin{equation}\label{eq1}
Tf=F \cdot(f\circ\varphi)\quad\quad(f\in H^p), \end{equation}
where $\varphi$ is a non-constant inner function (that is,
$|\varphi(z)|\leq1$ for $|z|\leq1$ and $|\varphi(z)|=1$ for $|z|=1$) and
$F\in H^p$ (we do not use any connection between $\varphi$ and $F$, but see
\cite[Theorem 1]{fo}). Moreover, $T\in\Iso(H^p)$ if and only if
\begin{equation}\label{eq2}
Tf=b\bigg(\frac{d\varphi}{dz}\bigg)^{1/p}(f\circ\varphi)\quad\quad(f\in
H^p), \end{equation}
where $b\in\mathbb T$ and $\varphi$ is a conformal map of the disc onto
itself.

\begin{proof} Let $T\in\Ref_{to}(\Iso(H^p))$. Since $T$ is an into
isometry,
there exists $F$ and $\varphi$ such that $T$ is of the form \eqref{eq1}
above. Observe that $F=T(1)$ and
$F\varphi=T(id)$ (here $id$ denotes the identity function on $U$). Take
$T_n, L_n\in\Iso(H^p)$ such that $F=\lim_nT_n(1)$ and
$F\varphi=\lim_nL_n(id)$. Taking into account the form of the
conformal maps of $U$ onto itself, we see that
\begin{eqnarray*}
F&=&\lim_{n\to\infty}b_n\bigg(\frac{1-|a_n|^2}{(1-\overline
a_nid)^2}\bigg)^{1/p}\\
F\varphi&=&\lim_{n\to\infty}\beta_n\bigg(\frac{1-|\alpha_n|^2}{(1-\overline
\alpha_nid)^2}\bigg)^{1/p}
\bigg(\frac{id-\alpha_n}{1-\overline\alpha_nid}\bigg),
\end{eqnarray*}
where $b_n,\beta_n\in\mathbb T$ and $a_n,\alpha_n\in U$ for all $n$.
Without loss of generality one can assume that the sequences $(b_n),
(\beta_n), (a_n), (\alpha_n)$ are convergent. Let
$$
b=\lim_{n\to\infty}b_n,\quad\beta=\lim_{n\to\infty}\beta_n,\quad
a=\lim_{n\to\infty}a_n,\quad\alpha=\lim_{n\to\infty}\alpha_n.
$$
Clearly, $|b|=|\beta|=1$. We show that $a,\alpha\in U$. Suppose that
$|a|=1$. Then, using the continuity of point-evaluations, we have
$$
F(z)=\lim_{n\to\infty}b_n\bigg(\frac{1-|a_n|^2}{(1-\overline
a_nz)^2}\bigg)^{1/p}=0
$$
for every $z\in U$. But this implies that
$T=0$, a contradiction. Therefore $a\in U$. The fact that $\alpha\in U$ can
be proved in a similar way: if we assume that $|\alpha|=1$, then
$F\varphi=0$ which means that $T$ is not injective. Consequently, we have
$\alpha \in U$. Now, it is easily seen that \begin{eqnarray*}
F(z)&=&b\bigg(\frac{1-|a|^2}{(1-\overline az)^2}\bigg)^{1/p}\\
(F\varphi)(z)&=&\beta
\bigg(\frac{1-|\alpha|^2}{(1-\overline \alpha z)^2}\bigg)^{1/p}
\bigg(\frac{z-\alpha}{1-\overline\alpha z}\bigg). \end{eqnarray*}
It remains to show that $a=\alpha$. Indeed, in this case we obtain that
$\varphi$ is a conformal map of $U$ onto itself, and that $T$ is of the
form \eqref{eq2} which gives us that $T$ is a surjective isometry.

Since $\varphi$ is an inner function on the disc, having in mind the form
of
$F$ and $F\varphi$, we infer that the expression $$
\frac{\beta}{b} \frac{z-\alpha}{1-{\overline{\alpha}}z} \bigg(
\frac{1-|\alpha |^2}{1-|a|^2}\bigg)^{1/p} \bigg( \frac{(1-\overline a
z)^2}{(1-\overline \alpha z)^2}\bigg)^{1/p} $$
is of modulus 1 whenever $|z|=1$.
It is now clear that the M\"{o}bius function $$
\omega(z)=\sqrt{\frac{1-|\alpha
|^2}{1-|a|^2}}\bigg(\frac{1-\overline a
z}{1-\overline \alpha z}\bigg)
$$
leaves $\mathbb T$ invariant. It is well-known that any holomorphic
function
on $U$ which is continuous on $\overline U$ and leaves $\mathbb T$
invariant
is either constant or has at least one zero in $U$. Since in case $a\neq
\alpha$ the only zero of $\omega$ is at $1/\overline a$ which is of modulus
greater than 1, it follows that $\omega$ is constant and hence we find that
$a=\alpha$. The proof is now complete. \end{proof}

\section{Banach algebras of holomorphic functions}\label{sbounded}
In this section we study the isometry group and the automorphism group of
Banach algebras of
holomorphic functions. If $\Omega \subset \mathbb C$ is a domain (that is,
an open, connected set), then $H(\Omega)$ denotes the algebra of all
holomorphic functions on $\Omega$ and $H^\infty(\Omega)$ stands for the
algebra of all bounded functions in $H(
\Omega)$. If $K \subset \mathbb C$ is a compact set, then $A(K)$
denotes the algebra of all continuous complex valued functions on
$K$ which are holomorphic
in the interior of $K$. On $H(\Omega)$ we consider the topology of the
uniform convergence on compact subsets, while $H^
\infty(\Omega)$ and $A(K)$ are equipped with the sup-norm topology.

Let $\cA$ be a unital semisimple commutative Banach algebra. $\cA$ is
called a uniform algebra if the spectral radius $r(.)$ is a complete norm
on $\cA$.
By the well-known fact on the uniqueness of Banach algebra norms on
semisimple Banach algebras \cite[6.1.1 Theorem]{Pal}, using the result that
an injective linear map between Banach spaces has closed range if and only
if it is bounded from below, it is easy to see that $\cA$ is a uniform
algebra if and only if the range of the Gelfand transformation on $\cA$ is
a
closed subalgebra of
the space of all continuous complex
valued functions on its structure space (see the
definition of sup-norm algebras in \cite{dLRW}). In what follows
$\sigma(.)$
denotes the spectrum.

\begin{theorem}\label{T:mol1}
Let $\cA$ be a uniform algebra. Every $T \in \Ref_{to}(\Iso(\cA))$ has the
form
\[
T(f)=\tau \psi(f) \qquad (f\in \cA)
\]
for some $\tau\in \cA$ with $\sigma(\tau)\subset \mathbb T$ and some unital
algebra endomorphism $\psi
:\cA \to \cA$. \end{theorem}

For the proof we recall the result
\cite[Theorem 3]{dLRW} stating that every surjective linear isometry of a
uniform algebra $\cA$ is an algebra automorphism of $\cA$ multiplied by an
element of $\cA$ whose spectrum is contained in $\mathbb T$.
\begin{proof}
Let $T\in \Ref_{to}(\Iso(\cA))$. Denote $\tau=T(1)$. Since \[
\sigma(\tau)=\{ \varphi(\tau) \, :\, \varphi \text{ is a character of } \cA
\},
\]
and the characters of $\cA$ are continuous,
we obtain $\sigma(\tau)\subset \mathbb T$. Consider the mapping
$\psi(.)=\tau^{-1}T(.)$. Let $\varphi$ be an arbitrary character of
$\cA$. We have $(\varphi \circ \psi)(1)=1$. Let $f\in \cA$ be an arbitrary
invertible element and choose $\tau_n \in \cA$ with $\sigma(\tau_n)\subset
\mathbb T$ and a sequence $(\psi_n)$ of automorphisms of $\cA$ such that
\[
\psi(f)=\lim_{n\to\infty} \tau^{-1}\tau_n \psi_n(f). \] We have
\begin{equation}\label{E:mol1}
\varphi(\psi(f))=\lim_{n\to\infty} \varphi(\tau)^{-1}\varphi(\tau_n)
\varphi(\psi_n(f)). \end{equation}
Clearly, the number $\varphi(\tau)^{-1}\varphi(\tau_n)$ is of modulus 1
and $\varphi(\psi_n(f))\in \sigma (f)$. Since $0\notin \sigma (f)$ and the
spectrum in a unital Banach algebra is compact, it follows that the limit
in
\eqref{E:mol1} is nonzero. Therefore, the linear functional $\varphi \circ
\psi$ maps 1 to 1 and it sends invertible elements to nonzero complex
numbers. By the well-known Gleason-Kahane-\. Zelazko theorem \cite[2.4.13
Theorem]{Pal} these imply that $\varphi \circ \psi$ is a multiplicative
linear functional for every character $\varphi$ of $\cA$. From the
semisimplicity of $\cA$ it follows that $\psi$ is an algebra
homomorphism.
\end{proof}

\begin{remark}\label{uniform}
In view of the above proof it is obvious that if $\cA$ is a (unital)
semisimple commutative Banach algebra, then every element of
$\Ref_{to}(\Aut( \cA ))$ is a (unital) algebra endomorphism of $\cA$.
\end{remark}

\begin{theorem}\label{T:mol2}
The isometry group and the automorphism group of the disc algebra
$A({\overline U})$ are topologically reflexive. \end{theorem}

Before the proof we recall some basic facts about the Banach algebra
$A({\overline U})$. First of all, the structure space of $A({\overline U})$
is ${\overline U}$ \cite[3.2.13]{Pal}. This gives us that $r(f)=\| f\|$ for
every $f\in A({\overline U})$ which shows that the disc algebra is a
uniform
algebra. Afterwards, one can easily verify that the automorphisms of
$A({\overline U})$ are precisely the maps of the form \[
f\longmapsto f\circ \varphi,
\]
where $\varphi: {\overline U} \to {\overline U}$ is a homeomorphism which
maps $U$ conformally onto itself.

\begin{proof}
Let $T \in \Ref_{to}(\Iso(A({\overline U})))$. By Theorem \ref{T:mol1}, $T$
is a unital endomorphism multiplied by a function $\tau\in A({\overline
U})$
whose spectrum is contained in $\mathbb T$. Since the spectrum of
$\tau$ is $\tau({\overline U})$, we find that $\tau$ must be a
constant function of modulus 1 and without loss of generality
we can assume
that $T$ is actually a unital endomorphism of the disc
algebra. It is not hard to verify that $T$ is necessarily of the form
\[
T(f)=f\circ \varphi \qquad (f \in A({\overline U})) \] for some
function $\varphi: {\overline U} \to {\overline U}$. Since
\[
T(id)\in \overline{\{ \lambda\psi(id) \, :\, |\lambda |=1, \psi \text{ is
an automorphism of the disc algebra}\}}, \]
taking into account the form of automorphisms of $A({\overline{U}})$ and
that
of the conformal maps of $U$ onto itself, it follows that there are
sequences $\lambda_n\in \mathbb T$ and $\alpha_n \in U$ such that
\[
\varphi(z)=T(id)(z)=
\lim_{n\to\infty} \lambda_n \frac{z-\alpha_n}{1-\overline{\alpha_n}z}, \]
where the convergence is uniform in $z\in {\overline U}$. Choosing
subsequences if necessary, we may assume that $(\lambda_n)$ and
$(\alpha_n)$
converge to $\lambda$ and $\alpha$, respectively. Suppose
that $\alpha$ is of modulus 1. Then \[
\lim_{n\to\infty} \frac{z-\alpha_n}{1-\overline{\alpha_n}z}=-\alpha \qquad
(|z| < 1). \]
Since the above convergence is uniform, there is an $n \in \mathbb N$ such
that
\[
\biggl| \frac{z-\alpha_n}{1-\overline{\alpha_n}z}+\alpha \biggr|
<\frac{1}{2} \qquad (|z| < 1).
\]
As the function $z \mapsto (z-\alpha_n)/(1-\overline{\alpha_n}z)$ maps $U$
onto itself, we arrive at a contradiction. Consequently, $\varphi$ is of
the
form
\[
\varphi(z)
= \lambda \frac{z-\alpha}{1-\overline{\alpha}z} \qquad (z\in
{\overline{U}})
\]
with $|\lambda|=1$ and $|\alpha|<1$, from which we obtain the surjectivity
of $T$. This completes the proof of the topological reflexivity of
$\Iso(A({\overline U}))$.

The statement about $\Aut(A({\overline
U}))$ follows from the same argument.
\end{proof}

\begin{theorem}\label{T:mol3}
Let $\Omega \subset \mathbb C$ be a simply connected domain. The isometry
group and the automorphism group of $H^\infty(\Omega)$ are topologically
reflexive.
\end{theorem}

\begin{proof} Clearly, the spectrum of any element $f$ of
$H^\infty (\Omega)$ is $\overline{f(\Omega)}$. Therefore, the spectral
radius is equal to the norm and hence $H^\infty (\Omega)$ is a uniform
algebra.

If $\Omega =\mathbb C$, then by Liouville's theorem we have
$H^\infty(\Omega)=\mathbb C$ and in this case the statement is trivial.

Suppose that $\Omega \subsetneq \mathbb C$. By the Riemann mapping theorem
we may assume that $\Omega = U$. First observe that since the structure
space of $H^\infty (U)$ contains $U$, it follows that every element of
$H^\infty (U)$ with spectrum contained in $\mathbb T$ is a constant
function
of modulus 1. Also, we know that every automorphism of the algebra
$H^\infty
(U)$ is induced by a conformal selfmap of $U$ \cite{Kak}. Now,
let $T \in \Ref_{to}(\Iso(H^\infty(U)))$.
Just as in the proof
of Theorem~\ref{T:mol2} we can assume that $T(1)=1$. One can check in a
way very similar to that we have followed there that $T(id)$ is a
conformal map of $U$ onto itself. Therefore, we can suppose that even
$T(id)=id$ holds true. With all these assumptions assume that
$z_0\in U$
and that $f\in H^\infty (U)$ is such that $f(z_0)=0$. We obtain \[
\frac{f}{id -z_0}\in H^\infty(U).
\]
Since $T$ is a unital algebra homomorphism, it follows that \[
\frac{T(f)}{id -z_0}=T\biggl(\frac{f}{id -z_0}\biggr)\in H^\infty(U) \]
which gives us that $T(f)(z_0)=0$. Let now $f\in H^\infty (U)$ be
arbitrary. Since the function $f-f(z_0)$ vanishes at $z_0$, it follows that
the same holds for $T(f-f(z_0))=T(f)-f(z_0)$. Thus, we have
$T(f)(z_0)=f(z_0)$. Since this is true for every $z_0 \in U$ and $f\in
H^\infty(U)$, we have $T(f)=f$ for all $f\in H^\infty(U)$. This completes
the proof of the theorem. \end{proof}

\begin{remark}
The main difference between the proofs of Theorem~\ref{T:mol2} and
Theorem~\ref{T:mol3} is that in the former one we were lucky to use the
form
of endomorphisms of $A({\overline{U}})$. Such an "inner" form for the
endomorphisms of $H^\infty (U)$ does not exist. The reason
is that the structure space of $H^\infty(U)$ is much bigger than
$U$.
\end{remark}

\begin{remark}\label{remsix} We have some remarks concerning the algebraic
reflexivity of
the isometry group and the automorphism group of algebras of holomorphic
functions on more general domains.

First recall that if $K\subset \mathbb C$ is a compact set whose complement
has finitely many components, then the structure space of $A(K)$ is just
$K$. In fact, this follows from the celebrated Mergelyan's approximation
theorem which states that in this case $A(K)$ is equal to the sup-norm
closure of the set of all rational functions with poles outside $K$
\cite[20.5 Theorem and Exercise 1, p. 427]{Rud}, \cite[3.2.14
Example]{Pal}.
We easily obtain that $A(K)$ is a uniform algebra and having a look at the
proofs of Theorem~\ref{T:mol1} and \ref{T:mol2}, one can easily verify that
the isometry group and the automorphism group of $A(K)$ are algebraically
reflexive.

Next, let
$\Omega\subset \mathbb C$ be a bounded domain. We assert that the isometry
group and the automorphism group of $H^\infty (\Omega)$ are algebraically
reflexive.
As for the proof, we can clearly assume
that
\[
\Omega \subset \{ z \in \mathbb C \, : \, \Re z>0\} \] and that $\Omega$
contains a real number. Obviously, if $\lambda \in \mathbb T$ is any number
different from 1, then there is a positive integer $n$ such that
$\lambda^n {\overline{\Omega}} \nsubseteq {\overline{\Omega}}$.
After this
short preparation, let $T: H^\infty(\Omega)\to H^\infty(\Omega)$ be a
local surjective isometry. Just as in the proof of Theorem~\ref{T:mol3}, we
can suppose that $T(1)=1$. Let $\lambda$ be a complex number of modulus 1
and let $\psi$ be an automorphism of $H^\infty (\Omega)$ such that
\begin{equation}\label{E:mol2}
T(id)=\lambda \psi(id).
\end{equation}
We note that it is not true for general domains that every automorphism of
$H^\infty(\Omega)$ is induced by a conformal selfmap of $\Omega$. We show
that in \eqref{E:mol2} we have $\lambda =1$. The spectrum of $id$ is
$\overline{\Omega}$. Clearly, every automorphism preserves the spectrum
while any unital algebra homomorphism is spectrum non-increasing. Hence,
from the relation above we infer that \[ \lambda{\overline{\Omega}}=
\lambda \sigma (id)=
\lambda \sigma(\psi(id))
=\sigma(T(id))\subset \sigma (id)
= \overline{\Omega}. \]
Since this implies that
$\lambda^n \overline{\Omega} \subset \overline{\Omega}$ for every $n\in
\mathbb N$, we obtain that $\lambda $ must equal 1. As the norm in
$H^\infty
(\Omega)$ is equal to the spectral radius, it follows that every
automorphism is a surjective isometry. Therefore, considering the map
$\psi^{-1} \circ T$ we can assume that our local surjective isometry
satisfies $T(1)=1$ and $T(id)=id$. Now, the proof can be completed as
in the last part of the proof of Theorem~\ref{T:mol3}. The algebraic
reflexivity of the automorphism group can be proved in an easier way.
\end{remark}

\section{Fr\'echet algebras of holomorphic functions}\label{sunbounded}
Consider now the full algebra $H(\Omega)$ equipped with the topology of
uniform convergence on compact subsets of $\Omega$. In this case
$H(\Omega)$, as a Fr\'echet space, is metrizable but there is no natural
metric on $H(\Omega)$ and, therefore, there is no natural notion of
isometry
for $H(\Omega)$. We have, however, the following result about the
automorphism group.
\begin{theorem}\label{T:mol4}
Let $\Omega\subset\mathbb C$ be a simply connected domain. Then the
automorphism group of $H(\Omega)$ is topologically reflexive if and only if
$\Omega\neq\mathbb C$.
\end{theorem}

\begin{proof}
The fact that the automorphism group of
$H(\mathbb C)$ is not topologically reflexive can be seen as follows.
Consider the automorphisms $\psi_n :H(\mathbb C) \to H(\mathbb C)$
defined by \[
(\psi_n(f))(z)=f(z/n) \qquad (z\in \mathbb C, n\in \mathbb N).
\]
Let $\psi(f)=f(0)$ $(f \in H(\mathbb C))$. Clearly, we have $\psi(f)=\lim_n
\psi_n(f)$ for every $f\in H(\mathbb C)$, but $\psi$ is not an automorphism
of $H(\mathbb C)$. This proves the ``only if part".

For the converse, recall that every every automorphism
of the algebra of all holomorphic functions on a domain is induced by a
conformal selfmap of the underlying set \cite{Ber, BS}. Just as above, we
can assume that $\Omega =U$.

Let $\psi\in \Ref_{to}(\Aut(H(U)))$. As in the previous section, it is easy
to see that any sequence of conformal mappings of $U$ onto itself has a
subsequence which converges (uniformly on compact subsets of $U$) to either
a map of the same kind or to a constant function of modulus 1. One can
readily check that this implies that for every $f\in A({\overline{U}})
\subset H(U)$ we have $\psi(f)\in A({\overline{U}})$ and the range of
$\psi(f)$ is contained in that of $f$. Therefore, by Gleason-Kahane-\.
Zelazko theorem $\psi$ is multiplicative on $A({\overline U})$ and just as
in the proof of Theorem~\ref{T:mol2} we infer that there is a function
$\varphi: {\overline{U}} \to {\overline{U}}$ such that \[
\psi(f)=f\circ \varphi \qquad ( f\in A({\overline{U}})). \]
Considering $\psi(id)$ we find that $\varphi$ is either a conformal map of
$U$ onto itself, or it is a constant function of modulus 1. Since
$\psi$ is supposed to be continuous and the polynomials are dense in
$H(U)$,
it follows that in the first case we have $\psi(f)=f\circ \varphi$ for
every
$f\in H(U)$, so in this case $\psi$ is an automorphism of $H(U)$. As for
the
second possibility, suppose that $\varphi \equiv \alpha$, where $\alpha $
is
a constant of modulus 1. Since the function $1/(\alpha -id)$ is a element
of
$H(U)$, by the continuity of $\psi$ we obtain \[
\psi \biggl(\frac{1}{\alpha-id}\biggr)=
\psi \biggl(\sum_n \frac{id^n}{\alpha^{n+1}}\biggr )= \sum_n
\frac{1}{\alpha}. \]
This contradiction shows that the second possibility cannot occur. The
proof is now complete.
\end{proof}

\begin{remark}
Let
$\Omega \subset \mathbb C$ be a bounded domain whose complement has
finitely
many components each of them with nonempty interior. We show that
the automorphism group of $H(\Omega)$ is algebraically reflexive. Let
$\psi:H(\Omega) \to H(\Omega)$ be a local automorphism. Since every
automorphism of $H(\Omega)$ is induced by a conformal selfmap of
$\Omega$, we see that
$\psi$ can be considered as a local automorphism of $H^\infty(\Omega)$. By
the second part of Remark~\ref{remsix}, it follows that $\psi$ is an
automorphism of $H^\infty(\Omega)$. Our topological condition on $\Omega$
was set to guarantee that, by Runge's approximation theorem \cite[13.9
Theorem]{Rud}, $H^\infty(\Omega)$ is dense in $H(\Omega)$. Therefore, by
the
continuity of $\psi$ we obtain that $\psi$ is a homomorphism of
$H(\Omega)$.
But the form of endomorphisms of $H(\Omega)$ is well-known. Using, for
example, \cite[Corollary]{BS}, we have $\psi(f)=f\circ \varphi$ $(f\in
H(\Omega))$. Since $\varphi=\psi(id)$ is, by assumption, a conformal map,
we
obtain that $\psi$ is an automorphism of $H(\Omega)$.

Finally, we
prove that the automorphism group of $H(\mathbb C)$ is also algebraically
reflexive. It is well-known that the conformal maps of $\mathbb C$ onto
itself are precisely the affine functions $z \mapsto az+b$, $a\neq 0$. Let
$\psi: H(\mathbb C) \to H(\mathbb C)$ be a local automorphism. Assume that
$\psi(id)=id$. By the continuity of $\psi$, it is sufficient to prove that
$\psi(id^n)=id^n$ holds for every $n\in \mathbb N$. Let $n\geq 3$. From the
equality $\psi(id +id^n)=\psi(id)+\psi(id^n)$ it follows that there are
complex numbers $a,b,c,d$ with $a,c\neq 0$ such that \[
(az+b)+(az+b)^n=z +(cz+d)^n \qquad (z\in \mathbb C). \] Applying an
appropriate affine transformation, the previous equality turns to \[
z+z^n=(a'z+b')+(c'z+d')^n,
\]
where $a',c'\neq 0$.
Comparing the coefficients of the polynomials above, since $1<n-1<n$,
it follows that $d'=0$ and then that $b'=0$. We obtain $a'=1, (c')^n=1$.
It is not hard to see that this yields $\psi(id^n)=id^n$. It remains to
check that $\psi(id^2)=id^2$. Let $\psi(id^2)=\varphi^2$, where
$\varphi$ is an affine function. Picking any $n\geq 4$ and using a very
similar argument as above but this time for
$\psi(id^2+id^n)=\psi(id^2)+\psi(id^n)$, we obtain \[ \varphi=
\frac{u_2}{u_n} id
\]
where $u_2$ is a second, while $u_n$ is an $n$th root of unity. Since this
holds for every $n\geq 4$, we readily have $\varphi^2=id^2$ and thus
$\psi(id^2)=id^2$. This completes the proof of the algebraic reflexivity of
$\Aut(H(\mathbb C))$.
\end{remark}

To conclude Sections 4, 5 and 6, we remark that it would be interesting to
extend our topological reflexivity results for more general domains on the
plane as well
as to give examples of exotic domains for which the corresponding spaces of
holomorphic functions have nonreflexive isometry groups, or automorphism
groups. It would be also interesting to treat similar problems for vector
valued holomorphic functions instead of scalar valued ones
(see \cite{Lin1, Lin2}).

\section{Spaces of continuous functions}\label{continuous}

In this section we give some examples concerning spaces of continuous
functions. It is proved in \cite[Theorem 2.2]{mz} that both $\Aut(C(K))$
and
$\Iso(C(K))$
are algebraically reflexive if $K$ is a first countable Hausdorff space.
Our
following example shows that reflexivity may fail even if $K$ lacks first
countability at only one point. Observe that if $\Gamma $ is a discrete
space, then for its one-point compactification $\Gamma \cup \{ \infty\}$ we
have
$C(\Gamma \cup \{ \infty\})=c(\Gamma)=\mathbb K \oplus c_0(\Gamma)$.

\begin{example}Let $\Gamma$ be an uncountable index set. Then
$\Iso(c(\Gamma))$ and $\Aut(c(\Gamma))$ are algebraically nonreflexive.
\end{example}

\begin{proof} First, the construction in Example \ref{uncountable} gives us
a non-surjective operator $T_0:c_0(\Gamma)\to c_0(\Gamma)$ which is both a
local surjective isometry and a local automorphism of $c_0(\Gamma)$. Since
$c(\Gamma)=\mathbb{K}\oplus c_0(\Gamma)$, we can extend $T_0$ to
$c(\Gamma)$
by defining
$$
T(\lambda1+f)=\lambda1+T_0(f)\quad(\lambda\in\mathbb K, f\in c_0(\Gamma)).
$$
Clearly, $T$ is a non-surjective local automorphism and hence a local
surjective isometry of $c(\Gamma)$. This completes the proof. \end{proof}

Our next example is a space of continuous functions whose set of local
automorphisms is as large as it can be according to Remark \ref{uniform}.
As
usual, we denote by $\mathbb N^*$ the growth of $\mathbb N$ in its Stone-\v
Cech compactification, that is, $\mathbb N^*=\beta \mathbb N \setminus
\mathbb N$ which is a compact Hausdorff space.

\begin{theorem} Every unital injective endomorphism of $C(\mathbb{N}^*)$ is
a local automorphism. Therefore, the isometry group and the automorphism
group of
$C(\mathbb{N}^*)$ fail to be algebraically reflexive.
\end{theorem}

\begin{proof} Let $\psi:C(\mathbb{N}^*)\to C(\mathbb{N}^*)$ be a unital
injective
endomorphism. As it is well-known, every unital endomorphism
of a $C(K)$ space is induced by a continuous selfmap of $K$. Hence there is
a surjective continuous mapping
$\varphi:\mathbb{N}^*\to\mathbb{N}^*$ for which $\psi (f)=f\circ\varphi$
for
all $f\in C(\mathbb{N}^*)$. Since $f$ and $\psi(f)$ have the same range,
there is an homeomorphism $\varphi_f$ of $\mathbb{N}^*$ such that
$\psi (f)=f\circ\varphi_f$ as a consequence of the following fact
about the structure of $\mathbb{N}^*$ \cite[p. 83]{Wa}: if $K$ is a compact
space of topological weight at most $\aleph_1$ and $f$ and $g$ are
continuous maps from $\mathbb{N}^*$ onto $K$, then there is a homeomorphism
$\phi$ of $\mathbb{N}^*$ such that $g=f\circ\phi$. Clearly, this yields
that
$\psi$
is both a local automorphism and a local surjective isometry, although
$\psi$
is never surjective unless $\varphi$ is injective. Since $\mathbb N^*$
contains two disjoint copies of itself which are clopen in $\mathbb N^*$
(see \cite[3.10, 3.14, 3.15]{Wa}), the existence of a
surjective non-injective mapping $\varphi: \mathbb{N}^*\to\mathbb{N}^*$ is
obvious. This completes the proof.
\end{proof}

Our final examples show that there is no much relation between the
reflexivity of the isometry group of $X$ and that of its dual $X^*$. This
solves in part Problem 2 at the end of \cite{mz}.

\begin{example}
$(a)$ There is a Banach space X such that $\Iso(X)$ is
topologically reflexive but $\Iso(X^*)$ fails to be algebraically
reflexive.
$(b)$ There is a Banach space Y with algebraically nonreflexive isometry
group such that $\Iso(Y^*)$ is algebraically reflexive.
\end{example}

\begin{proof} We first prove $(b)$. Let $Y=l_1(\Gamma)$, where $\Gamma$ is
an uncountable index set.
By Example~\ref{uncountable}, $\Iso(Y)$ is algebraically nonreflexive.
However, the isometry group of $Y^*=l_\infty(\Gamma)$ is algebraically
reflexive; this can be proved, {\it mutatis mutandis}, as in
the countable case (see \cite[Proposition 6]{bm}).

To verify $(a)$, let $K$ be Cook's continuum having no nontrivial
continuous
surjection onto itself \cite{Co}. We show that $\Iso(C(K))$ is topogically
reflexive. Let $T\in\Ref_{to}(\Iso(C(K)))$. By Theorem \ref{T:mol1}, one
has
$$Tf=\tau(f\circ\varphi),$$
where $\tau$ is unimodular and $\varphi$ is a continuous surjection of $K$
onto itself. Hence $\varphi$ is the
identity on $K$ and we obtain that $T$ is surjective. On the other hand, by
general representation results, one has $$
C(K)^*=l_1(K)\oplus_1L_1(\mu),
$$
where $\mu$ is a non-atomic measure. Since $K$ is uncountable (no continuum
is countable), it follows from Theorem \ref{lebesgue} that the isometry
group of $C(K)^*$ is algebraically nonreflexive. \end{proof}

We close the paper with the following open problem. Let $X$ be a locally
compact Hausdorff space and denote by $C_0(X)$ the Banach algebra of all
continuous complex valued functions on $X$ which vanish at infinity.
If $C_0(X)$ is separable (this means that the one-point compactification of
$X$ is metrizable), does it follow that $\Iso (C_0(X))$ is
algebraically reflexive? If this was the case,
we would get that the isometry group of any separable commutative
$C^*$-algebra is algebraically reflexive and we would have some hope to
obtain positive reflexivity results for the isometry group of some classes
of separable $C^*$-algebras.

\end{document}